\documentclass[11pt]{article}
\usepackage{amssymb}
\usepackage{amssymb}
\usepackage{indentfirst}
\usepackage{amscd}
\usepackage{amsmath}
\usepackage{latexsym}
\usepackage{amsfonts}
\usepackage{amssymb}
\usepackage{amsthm}
\usepackage{mathrsfs}
\usepackage{ams,amsmath,amsthm,amsfonts,amssymb,bm}

\begin{document}
\title{Operator inequalities dealing with operator equations
 \thanks
{ This work is supported by National Natural Science Fund of China
(10771011).}}
\author{ Jian Shi$^{1}$ \ \ Zongsheng Gao$^{2}$\hskip 1cm  \\
{\small $^{1,\ 2 }$LMIB $\&$ School of Mathematics and Systems Science,}\\
{\small Beihang University, Beijing, 100191, China}\\
}
         \date{}
         \maketitle

\maketitle \baselineskip 16pt \line(1,0){420}

\noindent{\bf Abstract}\ \ \ In this paper, we study the existence of  solutions of some kinds of
operator equations via operator inequalities. First, we
investigate characterizations of operator order $A\geqslant B
>0$ and chaotic operator order log $A \geqslant$ log $B$ for
positive definite operators $A$, $B$ in terms of operator equations, and generalize the results in
\cite{CSLin}. Then, we introduce applications of complete form of
Furuta inequality in operator equations. Some kinds of operator
equations are researched and related characterizations of solutions are proved. 
\vspace{0.2cm}


\vspace{0.2cm} \vskip   0.2cm \footnotetext[1]{Corresponding author,
E-mail addresses: shijian@ss.buaa.edu.cn (J. Shi).}   \footnotetext[2]{
E-mail addresses: zshgao@buaa.edu.cn (Z. Gao).}

\setlength{\baselineskip}{20pt}
%%%%%%%%%%%%%%%%%%%%%%%%%%%%%%%%%%%%%%%%%%%%%%%%%%first section %%%%%%%%%%%%%%%%%%%%%%%%%%%%%%%%%%%%%%%%%%%%
\section{Introduction }
\indent A bounded linear operator $A$ on a Hilbert space $H$ is said
to be positive if $\langle Ax,x \rangle \geqslant 0$ for all $x \in
H$. $A\geqslant 0$ and $A > 0$ mean a positive
semidefinite operator and a positive definite operator, respectively. If $X >0$, log $X$ is
defined as log $X= lim_{a\rightarrow +0} {\frac {X^{a}-I}{a}}$. We
call the relation log $A \geqslant$ log $B$ chaotic operator order
(in symbol:
$A\gg B$). \\
\indent In 1987 and 1995, Furuta showed the following operator
inequalities: \\
{\bf Theorem 1.1 }(Furuta Inequality)  \cite{TFuruta}. If $A\geqslant
B\geqslant O$, then for each $r\geqslant 0$,
$$
(B^{\frac r 2}A^{p}B^{\frac r 2})^{\frac 1 q} \geqslant (B^{\frac r
2}B^{p}B^{\frac r 2})^{\frac 1 q}, \eqno (1.1)
$$
$$
(A^{\frac r 2}A^{p}A^{\frac r 2})^{\frac 1 q} \geqslant (A^{\frac r
2}B^{p}A^{\frac r 2})^{\frac 1 q}  \eqno (1.2)
$$
hold for $p\geqslant 0$, $q\geqslant 1$ with $(1+r)q\geqslant
p+r$.\\
{\bf Theorem 1.2 }(Grand Furuta Inequality)  \cite{GFI}. If
$A\geqslant B\geqslant O$ with $A > O$, then for $t\in [0, 1]$ and
$p\geqslant 1$,
$$
A^{1-t+r}\geqslant (A^{\frac r 2}(A^{-{\frac t 2}}B^{p}A^{-{\frac t
2}})^{s}A^{\frac r 2})^{\frac {1-t+r}{(p-t)s+r}} \eqno (1.3)
$$
for $s\geqslant 1$ and $r\geqslant t$.\\
\indent In 1996, Tanahashi proved the conditions for $p$ and $q$ in
Furuta Inequality are the best possible if $r \geqslant 0$ in
\cite{BestFI}; In 2000, he proved the outer exponent value of Grand Furuta Inequality
 is the
best possible in \cite{BestGFI}. In 2008, Yuan et al.  proved a new operator
inequality, which is called Complete form of Furuta Inequality: \\
{\bf Theorem 1.3 }(Complete form of Furuta Inequality)  \cite{CFI}.
Let $A\geqslant B\geqslant O$, $r\geqslant 0$, $p>p_{0}\geqslant 0$,
and $s=min\{p, 2p_{0}+min\{1,r\}\}$. Then
$$
(A^{\frac r 2}B^{p_{0}}A^{\frac r 2})^{\frac
{s+r}{p_{0}+r}}\geqslant (A^{\frac r 2}B^{p}A^{\frac r 2})^{\frac
{s+r}{p+r}}.  \eqno (1.4)
$$
\indent Yuan et al.  also proved the optimality of the
outer exponent of Complete form of Furuta Inequality under the condition of $\{ 2p_{0}+min\{1,r\}\}\geq p$ or
$r\geq 1$ with $2p_{0}+1<p$ in \cite{CFI}. There were many related studies on Complete form of Furuta Inequality, such as \cite{OaM}, \cite{LAA}. \\
\indent Recently, some kinds of operator
equations have been shown and given deep discussion via Furuta inequality
and grand Furuta inequality (see \cite{CSLin}, \cite{JMAA}). In \cite{CSLin}, C.-S. Lin  proved the following two theorems, which the characterizations of operator order $A\geqslant B>0$
and chaotic operator order $A\gg B$ for positive
definite operator $A$, $B$ in terms of operator equations are discussed, respectively:\\
\noindent {\bf Theorem 1.4}  \cite{CSLin}. For $r\geqslant 0$ and a
nonnegative integer $n \geqslant 0$ such that $(1+r)(n+1)=p+r$ (so,
$p\geqslant
1 $), the following assertions are equivalent.\\
(A1) $A \geqslant B >O$ ;\\
(A2) $A^{1+r}\geqslant (A^{r/2}B^{p}A^{r/2})^{\frac {1}{n+1}}$ ;\\
(A3) $(B^{r/2}A^{p}B^{r/2})^{\frac {1}{n+1}}\geqslant B^{1+r}$ ;\\
(A4) there exists a unique $S > O$ with $\|  S \|
\leqslant 1$ such that\\
\indent $B^{p}=A^{1/2}S(A^{1+r}S)^{n}A^{1/2}=A^{1/2}(SA^{1+r})^{n}SA^{1/2}$ ;\\
(A5) there exists a unique $S > O$ with $\|S\|\leqslant  1$ such that\\
\indent $A^{p}=B^{1/2}S^{-1}(B^{1+r}S^{-1})^{n}B^{1/2}=B^{1/2}(S^{-1}B^{1+r})^{n}S^{-1}B^{1/2}$.\\
\noindent {\bf Theorem 1.5} \cite{CSLin}. For $p, r> 0$ and any
integer $n \geqslant 1$ such that $r(n+1)=p+r$, the following assertions are
equivalent.\\
(B1) $A\gg B$ ;\\
(B2) $A^{r}\geqslant (A^{r/2}B^{p}A^{r/2})^{\frac {1}{n+1}}$ ;\\
(B3) $(B^{r/2}A^{p}B^{r/2})^{\frac {1}{n+1}}\geqslant B^{r}$ ;\\
(B4) there exists a unique $T >O$ with $\| T\|
\leqslant 1$ such that\\
\indent $B^{p}=T(A^{r}T)^{n}=(TA^{r})^{n}T$ ;\\
(B5) there exists a unique $T >O$ with $\| T\|
\leqslant 1$ such that\\
\indent $A^{p}=T^{-1}(B^{r}T^{-1})^{n}=(T^{-1}B^{r})^{n}T^{-1}$.  \\ 
\\
\indent In this paper, we
will show generalized characterizations in terms of operator equations
 due to Theorem 1.4 and Theorem 1.5. Then, we will give
some applications of Complete form of Furuta Inequality in operator
equations. \\
\indent In order to prove our results, let's recall some well-known theorems:\\
\noindent {\bf Theorem 1.7} (L\"{o}wner-Heinz inequality)
\cite{Lowner}, \cite{Heinz}. $A\geqslant B\geqslant 0$ ensures
$A^{\alpha} \geqslant
B^{\alpha}$ for any $\alpha \in [0, 1]$. \\
\indent Theorem 1.7 assures that $A \geqslant B >0 \Rightarrow A\gg
B$.\\
\noindent {\bf Theorem 1.8} (Douglas's majorization and factorization theorem) \cite{Douglas}. The
following
assertions are equivalent for any operators $A$ and $B$.\\
\indent (a) range($A$)$\subseteq$ range ($B$);\\
\indent (b) $AA^{\ast }\leqslant \lambda BB^{\ast}$ for some $\lambda
\geqslant 0$ (majorization); \\
\indent (c) There exists a bounded operator C such that $A=BC$ (factorization).\\
Moreover, if (a), (b) and (c) are valid, then there exists a unique operator $C$ such that\\
\indent (i)$\|C\|^{2}=inf \{\mu |AA^{\ast }\leqslant BB^{\ast }\}$; \\
\indent (ii) Null($A$)=Null($C$);\\
\indent (iii) range($C$)$\subseteq $ range($B^{\ast}$).\\
\noindent {\bf Theorem 1.9} \cite{Yuan}. For $A, B > O$, $A \geqslant
B$ if and only if
$$
A^{1+t+r}\geqslant (A^{r/2}(A^{t/2}B^{p}A^{t/2})^{s}A^{r/2})^{
{(1+t+r)}/{((p+t)s+r)}}
$$
holds for $p\geqslant 1$, $t\geqslant 0$, $r\geqslant 0$ and
$s\geqslant {\frac {1+t}{p+t}}$.\\
\noindent {\bf Theorem 1.10} \cite{Yuan}. For $A, B > O$, log $A
\geqslant$ log  $ B$ if and only if
$$
A^{t+r}\geqslant (A^{r/2}(A^{t/2}B^{p}A^{t/2})^{s}A^{r/2})^{
{(t+r)}/{((p+t)s+r)}}
$$
holds for $p> 0$, $t\geqslant 0$, $r\geqslant 0$ and
$s\geqslant {\frac {t}{p+t}}$.\\

\section{The characterizations of operator order $A\geq B >O$ and chaotic operator order $A \gg B$}
\indent In this section , we will show further development of characterizations of operator order $A\geqslant B
>0$ and chaotic operator order log $A \geqslant$ log $B$ for
positive definite operators $A$, $B$ in terms of operator equations due to Lin's results in
\cite{CSLin}. \\
\indent  Throughout this section we assume that $A, B$ are positive
definite operators.\\
\noindent {\bf Theorem 2.1.} For $r\geqslant 0$, $t\geqslant 0$,
$p\geqslant 1$ and a nonnegative $n$ such that $(p+t)s+r=(n+1)\cdot
(1+t+r)$ (so, $s\geqslant {\frac {1+t}{p+t}}$), the following
assertions are
equivalent.\\
(C1) $A\geqslant B >O$ ; \\
(C2) $A^{1+t+r}\geqslant
(A^{r/2}(A^{t/2}B^{p}A^{t/2})^{s}A^{r/2})^{\frac {1}{n+1}}$ ;\\
(C3) $(B^{r/2}(B^{t/2}A^{p}B^{t/2})^{s}B^{r/2})^{\frac
{1}{n+1}}\geqslant B^{1+t+r}$ ;\\
(C4) there exists a unique $S > O$ with $\| S \|
\leqslant 1$ such that\\
\indent $B^{p}=A^{-{\frac t 2}}(A^{\frac
{1+t}{2}}(SA^{1+t+r})^{n}SA^{\frac {1+t}{2}})^{\frac 1 s}A^{-{\frac
t 2}} $\\
\indent \quad \ = $A^{-{\frac t 2}}(A^{\frac
{1+t}{2}}S(A^{1+t+r}S)^{n}A^{\frac
{1+t}{2}})^{\frac 1 s}A^{-{\frac t 2}}$ ;\\
(C5) there exists a unique $S > O$ with $\| S \|
\leqslant 1$ such that\\
\indent $A^{p}=B^{-{\frac t 2}}(B^{\frac
{1+t}{2}}(S^{-1}B^{1+t+r})^{n}S^{-1}B^{\frac {1+t}{2}})^{\frac 1
s}B^{-{\frac t 2}} $\\
\indent \quad \ = $B^{-{\frac t 2}}(B^{\frac
{1+t}{2}}S^{-1}(B^{1+t+r}S^{-1})^{n}B^{\frac
{1+t}{2}})^{\frac 1 s}B^{-{\frac t 2}}$.\\
\noindent {\bf Proof.} \\
(C1) $\Rightarrow $ (C2) is obvious by Theorem 1.9.\\
(C2) $\Rightarrow $ (C1). Take $n=0$, $t=0$, $r=0$ in
$(p+t)s+r=(n+1)\cdot(1+t+r)$ and (C2), then $ps=1$ and $A\geqslant
B$.\\
(C2) $\Leftrightarrow $ (C3) holds by $A\geqslant B >O
\Leftrightarrow B^{-1}\geqslant A^{-1} > O$
.\\
(C2) $\Rightarrow $ (C4). According to Theorem 1.8, there exists $C$
with $\| C\| \leqslant 1$ such that\\
$(A^{r/2}(A^{t/2}B^{p}A^{t/2})^{s}A^{r/2})^{\frac
{1}{2(n+1)}}=A^{\frac {1+t+r}{2}}C=C^{\ast }A^{\frac {1+t+r}{2}}$.\\
\indent Take $S=CC^{\ast }$, by $\| C\| \leqslant 1$ and $C$ is determined
by $A$ and $B$, we have $\| S \| \leqslant 1$, $S$ is unique, and
$(A^{r/2}(A^{t/2}B^{p}A^{t/2})^{s}A^{r/2})^{\frac {1}{n+1}}=A^{\frac
{1+t+r}{2}}SA^{\frac {1+t+r}{2}}$. Therefore, the following
equalities hold.
$$
 A^{r/2}(A^{t/2}B^{p}A^{t/2})^{s}A^{r/2} = (A^{\frac
{1+t+r}{2}}SA^{\frac {1+t+r}{2}})^{n+1}=A^{\frac
{1+t+r}{2}}(SA^{{1+t+r}})^{n}SA^{\frac {1+t+r}{2}} ;
 \eqno (2.1)
$$
$$
 A^{r/2}(A^{t/2}B^{p}A^{t/2})^{s}A^{r/2} = (A^{\frac
{1+t+r}{2}}SA^{\frac {1+t+r}{2}})^{n+1} =A^{\frac
{1+t+r}{2}}S(A^{{1+t+r}}S)^{n}A^{\frac {1+t+r}{2}}. \eqno (2.2)
$$
By (2.1) and (2.2), $(A^{t/2}B^{p}A^{t/2})^{s}=A^{\frac
{1+t}{2}}(SA^{{1+t+r}})^{n}SA^{\frac {1+t}{2}}
 =A^{\frac
{1+t}{2}}S(A^{{1+t+r}}S)^{n}A^{\frac {1+t}{2}}$, then (C4) is
obtained.\\
(C4) $\Rightarrow$ (C2).
\begin{eqnarray*}
& &(A^{r/2}(A^{t/2}B^{p}A^{t/2})^{s}A^{r/2})^{\frac {1}{n+1}}\\
&=&(A^{r/2}(A^{t/2}  A^{-{ t/ 2}}(A^{\frac
{1+t}{2}}(SA^{1+t+r})^{n}SA^{\frac {1+t}{2}})^{ 1/ s}A^{-{
t/ 2}}   A^{t/2})^{s}A^{r/2})^{\frac {1}{n+1}} \\
&=& ( A^{\frac {1+t+r}{2}}(SA^{{1+t+r}})^{n}SA^{\frac {1+t+r}{2}})^{\frac {1}{n+1}} \\
&=& ((A^{\frac {1+t+r}{2}}SA^{\frac {1+t+r}{2}})^{n+1})^{\frac
{1}{n+1}} \\
&=& A^{\frac {1+t+r}{2}}SA^{\frac {1+t+r}{2}} \\
&\leqslant &A^{ {1+t+r}}.
\end{eqnarray*}
The first equality is due to the equation of (C4), and the last
inequality is due to\\
 $O < S \leqslant \| S\|I \leqslant I$.\\
(C3) $\Rightarrow$ (C5). By (C3), we have $B^{-(1+t+r)}\geqslant
(B^{-r/2}(B^{-t/2}A^{-p}B^{-t/2})^{s}B^{-r/2})^{\frac {1}{n+1}}$.
According to Theorem 1.8, there exists $C$, such that $\|C\|\leqslant
1$, and $ (B^{-r/2}(B^{-t/2}A^{-p}B^{-t/2})^{s}B^{-r/2})^{\frac
{1}{2(n+1)}}= B^{-(1+t+r)/2}C=C^{\ast}B^{-(1+t+r)/2}$. Take
$S=CC^{\ast}$, because $\| C\| \leqslant 1$ and $C$ is uniquely
determined by $A$ and $B$, we have $\| S\| \leqslant 1$ and $S$ is
unique, and $(B^{-r/2}(B^{-t/2}A^{-p}B^{-t/2})^{s}B^{-r/2})^{\frac
{1}{n+1}}=B^{-(1+t+r)/2}SB^{-(1+t+r)/2} $. Therefore, the following
equality holds.
\begin{eqnarray*}
&\ &B^{r/2}(B^{t/2}A^{p}B^{t/2})^{s}B^{r/2} \\
&=&(B^{-r/2}(B^{-t/2}A^{-p}B^{-t/2})^{s}B^{-r/2})^{-1}\\
&=&(B^{-(1+t+r)/2}SB^{-(1+t+r)/2})^{-(n+1)}\\
&=&(B^{(1+t+r)/2}S^{-1}B^{(1+t+r)/2})^{n+1}\\
&=&B^{(1+t+r)/2}S^{-1}(B^{1+t+r}S^{-1})^{n}B^{(1+t+r)/2}\\
&=&B^{(1+t+r)/2}(S^{-1}B^{1+t+r})^{n}S^{-1}B^{(1+t+r)/2}.
\end{eqnarray*}
Then (C5) is obtained immediately.

(C5) $\Rightarrow $ (C3).
\begin{eqnarray*}
&\ &(B^{r/2}(B^{t/2}A^{p}B^{t/2})^{s}B^{r/2})^{\frac {1}{n+1}}\\
&=&(B^{r/2}(B^{t/2}  B^{-{ t/ 2}}(B^{\frac
{1+t}{2}}S^{-1}(B^{1+t+r}S^{-1})^{n}B^{\frac
{1+t}{2}})^{ 1/ s}B^{-{ t/ 2}}     B^{t/2})^{s}B^{r/2})^{\frac {1}{n+1}}\\
&=&(B^{\frac {1+t+r}{2}}S^{-1}(B^{ {1+t+r}}S^{-1})^{n}B^{\frac
{1+t+r}{2}})^{\frac {1}{n+1}}\\
&=&((B^{\frac {1+t+r}{2}}S^{-1}B^{\frac {1+t+r}{2}})^{n+1})^{\frac
{1}{n+1}}\\
&=&B^{\frac {1+t+r}{2}}S^{-1}B^{\frac {1+t+r}{2}}\\
&\geqslant &B^{ {1+t+r}}.
\end{eqnarray*}
The first equality is due to (C5), the last inequality is due to $O < S
\leqslant \| S \| I\leqslant I\Rightarrow S^{-1}
\geqslant I$.             \\
{\bf Remark 2.1.} The special case $r=0$ and $s=1$ of Theorem 2.1
just is Theorem 1.4. (A2), (A3), (A4), (A5) in Theorem 1.4 are extended to
(C2), (C3), (C4), (C5) in Theorem 2.1, respectively.\\
\noindent {\bf Theorem 2.2.} For $p$, $r > 0$, $t\geqslant 0$ and
any integer $n\geqslant 1$ such that $(p+t)s+r=(n+1)\cdot (t+r)$,
the following assertions are
equivalent.\\
(D1) $A\gg B$ ; \\
(D2) $A^{t+r}\geqslant
(A^{r/2}(A^{t/2}B^{p}A^{t/2})^{s}A^{r/2})^{\frac {1}{n+1}}$ ;\\
(D3) $(B^{r/2}(B^{t/2}A^{p}B^{t/2})^{s}B^{r/2})^{\frac
{1}{n+1}}\geqslant B^{t+r}$ ;\\
(D4) there exists a unique $T > O$ with $\| T \|
\leqslant 1$ such that\\
\indent $B^{p}=A^{-{\frac t 2}}(A^{\frac
{t}{2}}(TA^{t+r})^{n}TA^{\frac {t}{2}})^{\frac 1 s}A^{-{\frac t
2}}$\\
\indent \quad \ = $A^{-{\frac t 2}}(A^{\frac
{t}{2}}T(A^{t+r}T)^{n}A^{\frac
{t}{2}})^{\frac 1 s}A^{-{\frac t 2}}$ ;\\
(D5) there exists a unique $T > O$ with $\| T \|
\leqslant 1$ such that\\
\indent $A^{p}=B^{-{\frac t 2}}(B^{\frac
{t}{2}}(T^{-1}B^{t+r})^{n}T^{-1}B^{\frac {t}{2}})^{\frac 1
s}B^{-{\frac t 2}}$ \\
\indent \quad \ = $B^{-{\frac t 2}}(B^{\frac
{t}{2}}T^{-1}(B^{t+r}T^{-1})^{n}B^{\frac
{t}{2}})^{\frac 1 s}B^{-{\frac t 2}}$.\\
{\bf Proof.}\\
(D1) $\Rightarrow$ (D2) is obvious by Theorem 1.10.\\
(D2) $\Rightarrow$ (D1). Take $t=0$, $r={\frac 1 n}$, then $ps=1$
and $ A^{\frac 1 n} \geqslant (A^{\frac {1}{2n}}BA^{\frac
{1}{2n}})^{\frac {1}{n+1}}$. Take logarithm of both sides and
refining, we have  log $ A\geqslant {\frac {n}{n+1}} $ log
$(A^{\frac {1}{2n}}BA^{\frac {1}{2n}})$, and let
$n\rightarrow\infty$, log $A \geqslant $ log $B$ is obtained.\\
(D2) $\Leftrightarrow $ (D3) is due to the fact that $A\gg B
\Leftrightarrow B^{-1}\gg A^{-1}$.\\
The proofs of (D2) $\Leftrightarrow$ (D4) and (D3) $\Leftrightarrow
$ (D5) are similar to the proofs of (C2) $\Leftrightarrow $ (C4) and
(C3) $ \Leftrightarrow$ (C5), and we omit them here.      \\
{\bf Remark 2.2.} The special case $t=0$ and $s=1$ of Theorem 2.2
just is Theorem 1.5. (B2), (B3), (B4), (B5) in Theorem 1.5 are extended to
(D2), (D3), (D4), (D5) in Theorem 2.2, respectively.\\
{\bf Corollary 2.1.} For $p >0$, $t \geqslant 0$, the following
assertions are equivalent.\\
(E1) $A\gg B$ ; \\
(E2) $A^{p+t}\geqslant (A^{p/2}(A^{t/2}
B^{p}A^{t/2})^{(p+2t)/(p+t)}A^{p/2})^{1/2}$ ;\\
(E3) $(B^{p/2}(B^{t/2}
A^{p}B^{t/2})^{(p+2t)/(p+t)}B^{p/2})^{1/2}\geqslant B^{p+t}$ ;\\
(E4) there exists a unique $T > O$ with $\| T \|
\leqslant 1$ such that \\
\indent $B^{p}=A^{-t/2}(A^{t/2}TA^{p+t}TA^{t/2})^{(p+t)/(p+2t)}A^{-t/2}$ ;\\
(E5) there exists a unique $T > O$ with $\| T \|
\leqslant 1$ such that \\
\indent $A^{p}=B^{-t/2}(B^{t/2}T^{-1}B^{p+t}T^{-1}B^{t/2})^{(p+t)/(p+2t)}B^{-t/2}$.\\
{\bf Proof.} Take $n=1$, $r=p$ in Theorem 2.2, then $s=
(p+2t)/(p+t)$.                 \\
{\bf Remark 2.3.} If $t=0$, Corollary 2.1  just  is Corollary 3.1 in
\cite{CSLin}.\\
{\bf Corollary 2.2.} For any positive integer $m, n$, the following
assertions are equivalent.\\
(F1) $A\gg B$ ;\\
(F2) $A^{{\frac 1 m}+{\frac 1 n}}\geqslant (A^{\frac {1}
{2n}}(A^{\frac {1} {2m}}BA^{\frac {1} {2m}})^{\frac
{m+n+1}{m+1}}A^{\frac {1} {2n}})^{\frac {1}{n+1}}$ ;\\
(F3) $(B^{\frac {1} {2n}}(B^{\frac {1} {2m}}AB^{\frac {1}
{2m}})^{\frac
{m+n+1}{m+1}}B^{\frac {1} {2n}})^{\frac {1}{n+1}}\geqslant B^{{\frac 1 m}+{\frac 1 n}}$ ;\\
(F4) there exists a unique $T > O$ with $\| T \|
\leqslant 1$ such that \\
\indent $B=A^{-{\frac {1}{2m} }}(A^{\frac {1} {2m}}(TA^{{\frac 1
m}+{\frac 1 n}})^{n}TA^{\frac {1} {2m}})^{\frac
{m+1}{m+n+1}}A^{-{\frac {1}{2m} }}$ \\
\indent \quad =$A^{-{\frac {1}{2m} }}(A^{\frac {1} {2m}}T(A^{{\frac
1 m}+{\frac 1 n}}T)^{n}A^{\frac {1}
{2m}})^{\frac {m+1}{m+n+1}}A^{-{\frac {1}{2m} }}$;\\
(F5) there exists a unique $T > O$ with $\| T \|
\leqslant 1$ such that \\
\indent $A=B^{-{\frac {1}{2m} }}(B^{\frac {1} {2m}}(T^{-1}B^{{\frac
1 m}+{\frac 1 n}})^{n}T^{-1}B^{\frac {1} {2m}})^{\frac
{m+1}{m+n+1}}B^{-{\frac {1}{2m} }}$\\
\indent \quad $ = B^{-{\frac {1}{2m} }}(B^{\frac {1}
{2m}}T^{-1}(B^{{\frac 1 m}+{\frac 1 n}}T^{-1})^{n}B^{\frac {1}
{2m}})^{\frac {m+1}{m+n+1}}B^{-{\frac {1}{2m} }}$.\\
{\bf Proof.} Take $p=1$, $r={\frac 1 n}$, $t={\frac 1 m}$ in Theorem
2.2, then $s={\frac {m+n+1}{m+1}}$.      \\
{\bf Remark 2.4.} Let $m \rightarrow \infty$, Corollary 2.2 just is
Corollary 3.2 in \cite{CSLin}.\\
{\bf Corollary 2.3.} For any positive integer $m, n$ and some $a
>0$, the following
assertions are equivalent.\\
(G1) $a^{n+1}A\gg B$ ;\\
(G2) $a^{\frac {m+n+1}{m+1}}A^{{\frac 1 m}+{\frac 1 n}}\geqslant
(A^{\frac {1} {2n}}(A^{\frac {1} {2m}}BA^{\frac {1} {2m}})^{\frac
{m+n+1}{m+1}}A^{\frac {1} {2n}})^{\frac {1}{n+1}}$ ;\\
(G3) $a^{\frac {m+n+1}{m+1}}(B^{\frac {1} {2n}}(B^{\frac {1}
{2m}}AB^{\frac {1} {2m}})^{\frac
{m+n+1}{m+1}}B^{\frac {1} {2n}})^{\frac {1}{n+1}}\geqslant B^{{\frac 1 m}+{\frac 1 n}}$ ;\\
(G4) there exists a unique $T > O$ with $\|T \|
\leqslant 1$ such that \\
\indent $B=a^{n+1}A^{-{\frac {1}{2m} }}(A^{\frac {1}
{2m}}(TA^{{\frac 1 m}+{\frac 1 n}})^{n}TA^{\frac {1} {2m}})^{\frac
{m+1}{m+n+1}}A^{-{\frac {1}{2m} }}$ \\
\indent \quad =$a^{n+1}A^{-{\frac {1}{2m} }}(A^{\frac {1}
{2m}}T(A^{{\frac 1 m}+{\frac 1 n}}T)^{n}A^{\frac {1}
{2m}})^{\frac {m+1}{m+n+1}}A^{-{\frac {1}{2m} }}$ ;\\
(G5) there exists a unique $T > O$ with $\| T \|
\leqslant 1$ such that \\
\indent $a^{n+1}A=B^{-{\frac {1}{2m} }}(B^{\frac {1}
{2m}}(T^{-1}B^{{\frac 1 m}+{\frac 1 n}})^{n}T^{-1}B^{\frac {1}
{2m}})^{\frac {m+1}{m+n+1}}B^{-{\frac {1}{2m} }}$\\
\indent \quad \quad\quad=$B^{-{\frac {1}{2m} }}(B^{\frac {1}
{2m}}T^{-1}(B^{{\frac 1 m}+{\frac 1 n}}T^{-1})^{n}B^{\frac {1}
{2m}})^{\frac {m+1}{m+n+1}}B^{-{\frac {1}{2m} }}$.\\
{\bf Proof.} Replace $A$ by $a^{n+1}A$ in Corollary 2.2.       \\
{\bf Remark 2.5.} Let $m\rightarrow \infty$, Corollary 2.3 just is
Corollary 3.3 in \cite{CSLin}.\\
{\bf Corollary 2.4.} For any positive integer $m$ and some $a >0$,
the following
assertions are equivalent.\\
(H1) $a^{2}A\gg B$ ;\\
(H2) $a^{\frac {m+2}{m+1}}A^{1+{\frac 1 m}}\geqslant (A^{\frac {1}
{2}}(A^{\frac {1} {2m}}BA^{\frac {1} {2m}})^{\frac
{m+2}{m+1}}A^{\frac {1} {2}})^{\frac {1}{2}}$ ;\\
(H3) $a^{\frac {m+2}{m+1}}(B^{\frac {1} {2 }}(B^{\frac {1}
{2m}}AB^{\frac {1} {2m}})^{\frac
{m+2}{m+1}}B^{\frac {1} {2 }})^{\frac {1}{2}}\geqslant B^{1+{\frac 1 m}}$ ;\\
(H4) there exists a unique $T > O$ with $\| T \|
\leqslant 1$ such that \\
\indent $B=a^{2}A^{-{\frac {1}{2m} }}(A^{\frac {1} {2m}}TA^{1+{\frac
1 m}}TA^{\frac {1} {2m}})^{\frac
{m+1}{m+2}}A^{-{\frac {1}{2m} }}$ ;\\
(H5) there exists a unique $T > O$ with $\| T \|
\leqslant 1$ such that \\
\indent $a^{2}A=B^{-{\frac {1}{2m} }}(B^{\frac {1}
{2m}}T^{-1}B^{1+{\frac 1 m}}T^{-1}B^{\frac {1} {2m}})^{\frac
{m+1}{m+2}}B^{-{\frac {1}{2m} }}$ .\\
{\bf Proof.} Take $n=1$ in Corollary 2.3.            \\
{\bf Remark 2.6.} Let $m\rightarrow \infty$, Corollary 2.4 just is
Corollary
3.4 in \cite{CSLin}.\\
\section{Applications of Complete form of Furuta Inequality in Operator Equations}
\indent In this section,  we will introduce some applications of Complete form of
Furuta Inequality in operator equations. Several kinds of operator
equations will be researched and related characterizations of solutions will be proved.\\
\indent Throughout this section we assume $A$ and $B$ are positive
definite operators. \\
{\bf Theorem 3.1.} If $1\geqslant r\geqslant 0$, $p> p_{0}\geqslant
0$, for a nonnegative integer $n$ such that
$p+r=(n+1)((2p_{0}+r)+r)$, the following assertions are equivalent.\\
(3-1) $A\geqslant B >O$ ;\\
(3-2) $(A^{\frac r 2}B^{p_{0}}A^{\frac r 2})^{2}\geqslant (A^{\frac
r
2}B^{p}A^{\frac r 2})^{\frac {1}{n+1}}$ ;\\
(3-3) There exists a unique $S > O$ with $\| S\|
\leqslant 1$ such that \\
\indent $B^{p}=B^{p_{0}}A^{\frac r 2} S(A^{\frac r
2}B^{p_{0}}A^{r}B^{p_{0}}A^{\frac r 2}S)^{n}A^{\frac r 2}B^{p_{0}}
=B^{p_{0}}A^{\frac r 2} (SA^{\frac r
2}B^{p_{0}}A^{r}B^{p_{0}}A^{\frac r 2})^{n}SA^{\frac r
2}B^{p_{0}}$.\\
{\bf Proof.} (3-1) $\Rightarrow$ (3-2) is obvious by Theorem
1.3.\\
(3-2) $\Rightarrow$ (3-1). Take $p_{0}=0$, $r=1$, $n=0$ in
$p+r=(n+1)((2p_{0}+r+r)$ and (3-2), then $p=1$ and $A^{2}\geqslant
A^{\frac 1 2}BA^{\frac 1 2}$.\\
(3-2) $\Rightarrow$ (3-3). According to Theorem 1.8, there exists an
operator $C$ with $\| C\| \leqslant 1 $ such that $(A^{\frac r
2}B^{p}A^{\frac r 2})^{\frac {1}{2(n+1)}}= (A^{\frac r
2}B^{p_{0}}A^{\frac r 2})C=C^{\ast }(A^{\frac r 2}B^{p_{0}}A^{\frac
r 2})$. Take $S=CC^{\ast}$, then $\| S\| =\| CC^{\ast} \| \leqslant
1$ and
$$
(A^{\frac r 2}B^{p}A^{\frac r 2})^{\frac {1}{n+1}} = (A^{\frac r
2}B^{p_{0}}A^{\frac r 2}) S (A^{\frac r 2}B^{p_{0}}A^{\frac r 2}).
\eqno (3.1)
$$
\indent By (3.1), $S$ is unique and the following equality holds,
\begin{eqnarray*}
&\ &A^{\frac r 2}B^{p}A^{\frac r 2}\\
&=& ((A^{\frac r 2}B^{p_{0}}A^{\frac
r 2}) S (A^{\frac r 2}B^{p_{0}}A^{\frac r 2}))^{n+1}\\
&=& A^{\frac r 2}B^{p_{0}}A^{\frac r 2} S((A^{\frac r
2}B^{p_{0}}A^{\frac r 2})^{2}S)^{n}A^{\frac r 2}B^{p_{0}}A^{\frac r
2}\\  &=& A^{\frac r 2}B^{p_{0}}A^{\frac r 2} (S(A^{\frac r
2}B^{p_{0}}A^{\frac r 2})^{2})^{n}SA^{\frac r 2}B^{p_{0}}A^{\frac r
2}.
\end{eqnarray*}
\indent By above equalities, we can obtain (3.3).\\
(3-3) $\Rightarrow$ (3-2).
\begin{eqnarray*}
&\ &(A^{\frac r 2}B^{p}A^{\frac r 2})^{\frac {1}{n+1}}\\ &=&
(A^{\frac r 2} B^{p_{0}}A^{\frac r 2} S(A^{\frac r
2}B^{p_{0}}A^{r}B^{p_{0}}A^{\frac r 2}S)^{n}A^{\frac r 2}B^{p_{0}}A^{\frac r 2})^{\frac {1}{n+1}} \\
&=& ((A^{\frac r 2} B^{p_{0}}A^{\frac r 2}SA^{\frac r 2}
B^{p_{0}}A^{\frac r 2})^{n+1})^{\frac {1}{n+1}} \\
&=& (A^{\frac r 2} B^{p_{0}}A^{\frac r 2})S(A^{\frac r 2}
B^{p_{0}}A^{\frac r 2})\\
&\leqslant  & (A^{\frac r 2} B^{p_{0}}A^{\frac r 2})^{2}.
\end{eqnarray*}
\indent The first equality is due to (3-3), and the inequality is
due
to $O< S \leqslant \| S\| I\leqslant I$.  \\
\\
{\bf Corollary 3.1.} If $1\geqslant r\geqslant 0$, $p>
p_{0}\geqslant 0$, for a nonnegative integer $n \geqslant 0$ such that
$p+r=(n+1)((2p_{0}+r)+r)$, the following assertions are equivalent.\\
(3-1) $A\geqslant B >O$ ;\\
(3-4) $(B^{\frac r 2}A^{p}B^{\frac r 2})^{\frac {1}{n+1}}\geqslant
(B^{\frac r
2}A^{p_{0}}B^{\frac r 2})^{2}$ ;\\
(3-5) There exists a unique $S > O$ with $\| S\|
\leqslant 1$ such that \\
$A^{p}=A^{p_{0}}B^{\frac r 2} S^{-1}(B^{\frac r
2}A^{p_{0}}B^{r}A^{p_{0}}B^{\frac r 2}S^{-1})^{n}B^{\frac r
2}A^{p_{0}} =A^{p_{0}}B^{\frac r 2} (S^{-1}B^{\frac r
2}A^{p_{0}}B^{r}A^{p_{0}}B^{\frac r 2})^{n}S^{-1}B^{\frac r
2}A^{p_{0}}$.\\
{\bf Proof.} Replace $A$ and $B$ by $B^{-1}$ and $A^{-1}$ ,
respectively, in (3-2) and (3-3), then  (3-4) and (3-5) are obtained.   \\ \\
{\bf Theorem 3.2.}  If $r \geqslant 1$, $p> p_{0}\geqslant 0$, for a
nonnegative integer $n$ such that
$p+r=(n+1)((2p_{0}+1)+r)$, the following assertions are equivalent.\\
(3-1) $A\geqslant B >O$ ;\\
(3-6) $(A^{\frac r 2}B^{p_{0}}A^{\frac r 2})^{\frac
{2p_{0}+1+r}{p_{0}+r}}\geqslant (A^{\frac r
2}B^{p}A^{\frac r 2})^{\frac {1}{n+1}}$ ;\\
(3-7) There exists a unique $S > O$ with $\| S\| \leqslant 1$ such
that
\begin{eqnarray*} B^{p} &=& A^{-{\frac r 2}}(A^{\frac r 2}B^{p_{0}}A^{\frac r
2})^{\frac {2p_{0}+1+r}{2p_{0}+2r}} S((A^{\frac r
2}B^{p_{0}}A^{\frac r 2})^{\frac {2p_{0}+1+r}{ p_{0}+
r}}S)^{n}(A^{\frac r 2}B^{p_{0}}A^{\frac r 2})^{\frac
{2p_{0}+1+r}{2p_{0}+2r}} A^{-{\frac r 2}} \\ &=& A^{-{\frac r
2}}(A^{\frac r 2}B^{p_{0}}A^{\frac r 2})^{\frac
{2p_{0}+1+r}{2p_{0}+2r}} (S(A^{\frac r 2}B^{p_{0}}A^{\frac r
2})^{\frac {2p_{0}+1+r}{p_{0}+r}})^{n}S(A^{\frac r
2}B^{p_{0}}A^{\frac r 2})^{\frac {2p_{0}+1+r}{2p_{0}+2r}}
 A^{-{\frac r 2}}.
\end{eqnarray*}
{\bf Proof.} (3-1) $\Rightarrow$ (3-6) is obvious by Theorem 1.3.\\
(3-6) $\Rightarrow$ (3-1). Take $n=0$, $p_{0}=0$, $r=1$ in
$p+r=(n+1)((2p_{0}+1)+r)$ and (3-6), then $p=1$ and $A^{2} \geqslant
A^{\frac 1 2}B A^{\frac 1 2}$. (3-1) is obtained.\\
\indent The proof of (3-6) $\Leftrightarrow$ (3-7) is similar to the
proof of
(3-2) $\Leftrightarrow$ (3-3) in Theorem 3.1, so we omit it here.  \\ \\
{\bf Corollary 3.2.}  If $r \geqslant 1$, $p> p_{0}\geqslant 0$, for a
nonnegative integer $n$ such that
$p+r=(n+1)((2p_{0}+1)+r)$, the following assertions are equivalent.\\
(3-1) $A\geqslant B >O$ ;\\
(3-8) $(B^{\frac r 2}A^{p}B^{\frac r 2})^{\frac {1}{n+1}}\geqslant
(B^{\frac r 2}A^{p_{0}}B^{\frac r 2})^{\frac
{2p_{0}+1+r}{p_{0}+r}}$ ;\\
(3-9) There exists a unique $S > O$ with $\| S\| \leqslant 1$ such
that
\begin{eqnarray*} A^{p} &=& B^{-{\frac r 2}}(B^{\frac r 2}A^{p_{0}}B^{\frac r
2})^{\frac {2p_{0}+1+r}{2p_{0}+2r}} S^{-1}((B^{\frac r
2}A^{p_{0}}B^{\frac r 2})^{\frac
{2p_{0}+1+r}{p_{0}+r}}S^{-1})^{n}(B^{\frac r 2}A^{p_{0}}B^{\frac r
2})^{\frac {2p_{0}+1+r}{2p_{0}+2r}} B^{-{\frac r 2}} \\ &=&
B^{-{\frac r 2}}(B^{\frac r 2}A^{p_{0}}B^{\frac r 2})^{\frac
{2p_{0}+1+r}{2p_{0}+2r}} (S^{-1}(B^{\frac r 2}A^{p_{0}}B^{\frac r
2})^{\frac {2p_{0}+1+r}{p_{0}+r}})^{n}S^{-1}(B^{\frac r
2}A^{p_{0}}B^{\frac r 2})^{\frac {2p_{0}+1+r}{2p_{0}+2r}} B^{-{\frac
r 2}}.
\end{eqnarray*}
{\bf Proof.} Replaced $A$ and $B$ by $B^{-1}$ and $A^{-1}$ in (3-6)
and (3-7) of Theorem 3.2, respectively, then (3-8) and (3-9) are obtained.  \\ \\
{\bf Theorem 3.3.} If $r\geqslant 0$, $2p_{0} + min\{1, r\}\geqslant
p >p_{0}\geqslant 0$, for a positive integer $n$ such that
$n(p+r)=(n+1)(p_{0}+r)$, the following assertions are equivalent.\\
(3-1) $A\geqslant B >O$ ;\\
(3-10) $(A^{\frac r 2}B^{p_{0}}A^{\frac r 2})^{1+{\frac 1
n}}\geqslant
A^{\frac r 2}B^{p}A^{\frac r 2} $ ;\\
(3-11) There exists a unique $S>O$ with $\| S\|
\leqslant 1$ such that \\
\indent $(A^{\frac r 2}B^{p_{0}}A^{\frac r 2})^{n+1}= ((A^{\frac r
2}B^{p}A^{\frac r 2})^{\frac 1 2}S^{-1}(A^{\frac r 2}
B^{p}A^{\frac r 2})^{\frac 1 2})^{n}$ .   \\
{\bf Proof.} (3-1) $\Rightarrow$ (3-10) is obvious by Theorem 1.3.\\
(3-10) $\Rightarrow$ (3-1). Take $n=1$, $p_{0}=0$, $r=1$ in
$n(p+r)=(n+1)(p_{0}+r)$ and (3-10), then $p=1$ and $A^{2} \geqslant
A^{\frac 1 2}B A^{\frac 1 2}$. $A \geqslant B$ is obtained.\\
(3-10) $\Rightarrow$ (3-11). By (3-10) we have $A^{-{\frac r
2}}B^{-p}A^{-{\frac r 2}} \geqslant (A^{-{\frac r
2}}B^{-p_{0}}A^{-{\frac r 2}})^{\frac {n+1}{n}}$. According to
Theorem 1.8, there exists an operator $C$ with $\| C\|
\leqslant 1$, such that \\
\indent $(A^{-{\frac r 2}}B^{-p_{0}}A^{-{\frac r 2}})^{\frac
{n+1}{2n}} =(A^{-{\frac r 2}}B^{-p}A^{-{\frac r 2}})^{\frac 1
2}C=C^{\ast}(A^{-{\frac r 2}}B^{-p}A^{-{\frac r 2}})^{\frac 1 2}$.
\\
Take $S=CC^{\ast}$, then $S > O$,  $\| S\| \leqslant 1$ and the
following equality,
$$
(A^{-{\frac r 2}}B^{-p_{0}}A^{-{\frac r 2}})^{\frac {n+1}{n}}
=(A^{-{\frac r 2}}B^{-p}A^{-{\frac r 2}})^{\frac 1 2}S(A^{-{\frac r
2}}B^{-p}A^{-{\frac r 2}})^{\frac 1 2}. \eqno (3.2)
$$
By (3.2), $S$ is unique and $(A^{\frac r 2}B^{p_{0}}A^{\frac r
2})^{n+1}= ((A^{\frac r 2}B^{p}A^{\frac r 2})^{\frac 1
2}S^{-1}(A^{\frac r 2} B^{p}A^{\frac r 2})^{\frac 1 2})^{n}$.\\
(3-11) $\Rightarrow$ (3-10). Because of (3-11) and the fact that $S
>O$ with $\| S\| \leqslant 1 \Rightarrow S^{-1} \geqslant I$, we
have $(A^{\frac r 2}B^{p_{0}}A^{\frac r 2})^{1+{\frac 1 n}}
=(A^{\frac r 2}B^{p}A^{\frac r 2})^{\frac 1 2}S^{-1}(A^{\frac r 2}
B^{p}A^{\frac r 2})^{\frac 1 2}
 \geqslant  A^{\frac r 2} B^{p}A^{\frac r 2}$.   \\ \\
{\bf Corollary 3.3.} If $r\geqslant 0$, $2p_{0} + min\{1,
r\}\geqslant p
>p_{0}\geqslant 0$, for a positive integer $n$ such that
$n(p+r)=(n+1)(p_{0}+r)$, the following assertions are equivalent.\\
(3-1) $A\geqslant B >O$ ;\\
(3-12) $B^{\frac r 2}A^{p}B^{\frac r 2} \geqslant
(B^{\frac r 2}A^{p_{0}}B^{\frac r 2})^{1+{\frac 1 n}}$ ;\\
(3-13) There exists a unique $S>O$ with $\| S\|
\leqslant 1$ such that \\
\indent $(B^{\frac r 2}A^{p_{0}}B^{\frac r 2})^{n+1}= ((B^{\frac r
2}A^{p}B^{\frac r 2})^{\frac 1 2}S(B^{\frac r 2}
A^{p}B^{\frac r 2})^{\frac 1 2})^{n}$ .\\
{\bf Proof.} Replace $A$ and $B$ by $B^{-1}$ and $A^{-1}$ ,
respectively, in (3-10) and (3-11) of Theorem 3.3, then  (3-12) and (3-13) are
obtained.\\

\begin{center}

\end{center}
\end{document}